\documentclass[11pt]{article} 
\usepackage[leqno]{amsmath}
\usepackage{amssymb,amsthm,url}
\usepackage[letterpaper,margin=1in]{geometry}
\usepackage{graphicx}
\usepackage{color}
\numberwithin{equation}{section}
\numberwithin{figure}{section}
\usepackage{enumitem}		

\usepackage[style=authoryear,backend=bibtex,maxbibnames=99]{biblatex}
\newlength{\lyxlabelwidth}      

\definecolor{GREEN}{RGB}{0, 180, 0}
\definecolor{BLUE}{RGB}{0, 0, 180}

\DeclareMathOperator{\Proj}{Proj}

\DeclareMathOperator{\Des}{Des}

\DeclareMathOperator{\Peak}{Peak}

\DeclareMathOperator{\desc}{desc}
\DeclareMathOperator{\anc}{anc}

\newcommand\cpp{\mathbf{T}}

\newcommand\f{\mathbf{f}}
\newcommand\fo{\overset{\approx}{\f}}

\newcommand\calh{\mathcal{H}}
\newcommand\calhn{\mathcal{H}_n}

\newcommand\barcalh{\bar{\mathcal{H}}}

\newcommand\calb{\mathcal{B}}
\newcommand\calbn{\mathcal{B}_n}
\newcommand\hatcalb{\check{\mathcal{B}}}
\newcommand\hatcalbn{\check{\mathcal{B}}_n}

\newcommand\barcalb{\bar{\mathcal{B}}}

\newcommand\proda{m^{[a]}}
\newcommand\coproda{\Delta^{[a]}}

\newcommand\calp{\mathcal{P}}

\newcommand\sn{\mathfrak{S}_n}

\newcommand\sk{\mathfrak{S}_k}

\newcommand\sj{\mathfrak{S}_j}

\date{}

\usepackage{enumitem}		
\theoremstyle{plain}
\newtheorem{thm}{\protect\theoremname}[section]
  \theoremstyle{definition}
  \newtheorem{defn}[thm]{\protect\definitionname}
  \theoremstyle{plain}
  \newtheorem{prop}[thm]{\protect\propositionname}
  \theoremstyle{plain}
  
  \theoremstyle{plain}
  
  \theoremstyle{remark}
  \newtheorem*{rem*}{\protect\remarkname}
  \theoremstyle{definition}
  \newtheorem{example}[thm]{\protect\examplename}
\newenvironment{elabeling}[2][]%
{\settowidth{\lyxlabelwidth}{#2}
\begin{description}[font=\normalfont,style=sameline,
leftmargin=\lyxlabelwidth,#1]}
{\end{description}}

\makeatother

  \providecommand{\definitionname}{Definition}
  \providecommand{\examplename}{Example}
  \providecommand{\propositionname}{Proposition}
\providecommand{\theoremname}{Theorem}

\begin{document}

\title{Card-Shuffling via Convolutions of Projections on Combinatorial Hopf
Algebras}

\author{C. Y. Amy Pang}

\maketitle

\begin{abstract}
Recently, Diaconis, Ram and I created Markov chains out of the coproduct-then-product
operator on combinatorial Hopf algebras. These chains model the breaking
and recombining of combinatorial objects. Our motivating example was
the riffle-shuffling of a deck of cards, for which this Hopf algebra
connection allowed explicit computation of all the eigenfunctions.
The present note replaces in this construction the coproduct-then-product
map with convolutions of projections to the graded subspaces, effectively
allowing us to dictate the distribution of sizes of the pieces in
the breaking step of the previous chains. An important example is
removing one ``vertex'' and reattaching it, in analogy with top-to-random
shuffling. This larger family of Markov chains all admit analysis
by Hopf-algebraic techniques. There are simple combinatorial expressions
for their stationary distributions and for their eigenvalues and multiplicities
and, in some cases, the eigenfunctions are also calculable.

$ $

R\'{e}cemment, avec Diaconis et Ram, nous avons construit des chaines de Markov \`{a} partir de l'op\'{e}rateur ``coproduit-puis-produit"  d\'{e}fini sur un alg\`{e}bre de Hopf combinatoire. Ces chaines mod\'{e}lisent la d\'{e}construction et la construction d'objets combinatoires. La motivation \'{e}tait le ``m\'{e}lange \`{a} l'am\'{e}ricaine", une m\'{e}thode populaire pour m\'{e}langer un jeu de cartes, pour lequel les liens avec les alg\`{e}bres de Hopf combinatoires nous a permis de calculer explicitement toutes les fonctions propres. Ici, on g\'{e}n\'{e}ralise cette construction en rempla\c{c}ant l'op\'{e}rateur ``coproduit-puis-produit" par les convolutions de projections sur les composantes gradu\'{e}es de l'alg\`{e}bre. Ceci nous permet de stipuler les tailles des pi\`{e}ces dans la d\'{e}composition des objets combinatoires. Un exemple important est la suppression et l'insertion d'un ``sommet'', par analogie avec la biblioth\`{e}que de Tsetlin. On constate que toutes ces chaines peuvent \^{e}tre analys\'{e}es par des techniques provenant de la th\'{e}orie des alg\`{e}bres de Hopf combinatoires. On prouve des expressions combinatoires simples pour les distributions stationnaires ainsi que pour les valeurs propres et leurs multiplicit\'{e}s. Dans certains cas, il est possible de calculer les fonctions propres associ\'{e}es.
\end{abstract}

This version: \today. This extended abstract was accepted as a talk
for the 27th International Conference on Formal Power Series and Algebraic
Combinatorics (FPSAC) in Daejoen, South Korea, in July 2015, and is
published in a proceedings volume of Discrete Mathematics and Theoretical
Computer Science (DMTCS).

\section{Introduction\label{sec:intro}.}

\subsection*{Background: Markov chains from Hopf-powers}

Possibly the most popular model of card-shuffling is the Gilbert-Shannon-Reeds
(GSR) \emph{riffle-shuffle}: cut the deck into two piles according
to the (symmetric) binomial distribution, then drop one-by-one the
bottom card from either pile, chosen with probability proportional
to the current pile size. This second step is equivalent to all interleavings
of the two piles (counted with multiplicity) being equally likely. 

Amongst the plethora of results concerning this shuffle, the most
notable must be the work of \parencite{originalriffleshuffle}, who determined
that $\frac{3}{2}\log n$ shuffles are necessary and sufficient to
randomise a deck of $n$ cards. Central to their argument is the generalisation
of the GSR model to \emph{$a$-handed shuffles}: cut the deck into
\emph{a} piles according to a (symmetric) multinomial distribution,
then drop the bottom cards from a pile chosen with probability proportional
to pile size as before. Performing the GSR shuffle $t$ times is then
the same as a $2^{t}$-handed shuffle, so analysing long-term behaviour
of 2-handed riffle-shuffles is equivalent to letting the number of
hands tend to infinity.

 \parencite{hopfpowerchains} observed that the transition probabilities
of the riffle-shuffle are, up to scaling, the coefficients of the
coproduct-then-product map $m\Delta$ on the shuffle algebra. Furthermore,
the coefficients of the \emph{$a$th Hopf-power map} $\proda\coproda$
on the shuffle algebra give the transition probabilities of the $a$-handed
riffle-shuffle. We then defined Markov chains on the bases of other
combinatorial Hopf algebras by setting their transition probabilities
to be such Hopf-power coefficients, with a little modificiation (via
the Doob transform). These \emph{Hopf-power Markov chains} model the
breaking then recombining of the combinatorial objects indexing the
bases of the algebras. The thesis of \parencite{mythesis} greatly extends
the Hopf-power Markov chain framework; this encompasses a restriction-then-induction
chain on  representations of the symmetric groups, and a tree-pruning
model - see Examples \ref{ex:repsofsn} and \ref{ex:trees} below.

\subsection*{Probabilistic conclusions from Hopf-algebraic techniques}

The benefit of this viewpoint is two-fold. Firstly, as recorded by
\parencite[Th. 4.7.1]{mythesis}, maps between Hopf algebras which ``respect
the bases'' induce projections of the related chains. \parencite{hpmccompositions}
applied this to a map from the shuffle algebra to the algebra of quasisymmetric
functions, to conclude that the positions of descents under riffle-shuffling
of a deck of distinct cards is a Markov statistic. (A \emph{descent}
occurs where a card has greater value than the card immediately below
it.) This means that the probability of a shuffle producing descents
in prescribed positions depends only on the positions of descents
before the shuffle, and not on the exact deck order. (This fact also
follows from the descent set being a ``shuffle-compatible statistic'',
which \parencite{gesseltalk} attributes to Stanley.) \parencite[Sec. 4.7]{mythesis}
constructs many Markov statistics for inverse riffle-shuffling out
of  commutation quotients of the free associative algebra.

The second way in which the Hopf formulation aids in studying these
Markov chains is that, in many cases, there are algorithms to compute
a basis of eigenvectors for the Hopf-power maps and hence the transition
matrices. This gives interesting information about the long term behaviour
of the chain. As an example, \parencite[Prop. 6.1.3 and Prop. 6.1.5]{mythesis}
state that, if a deck of $n$ distinct cards was originally in ascending
order, then, after $t$ iterations of the $a$-handed shuffle, the
expected number of descents is $(1-a^{-t})\frac{n-1}{2}$, and the
expected number of peaks is $(1-a^{-2t})\frac{n-2}{3}$. (A \emph{peak}
is a triple of adjacent cards with the middle one having greatest
value.) Although the algorithms do not provide all eigenvectors of
all Hopf-power Markov chains, their stationary distributions are always
computable.

\subsection*{A new extension: Markov chains from convolutions of projections}

As \parencite[Th. 4.4.1]{mythesis} shows, the breaking step of a Hopf-power
Markov chain always involves a symmetric multinomial distribution.
However, it is sometimes more natural to consider other distributions.
For example, the restriction-then-induction chain mentioned above
is difficult to express in terms of partitions, because the restriction
of a symmetric group representation to a multinomially-chosen Young
subgroup involves Littlewood-Richardson coefficients. Simpler is the
Markov chain which removes a random corner box and re-inserts it in
a random position.

The discovery in this extended abstract is that, by replacing $\proda\coproda$
in the definition of a Hopf-power Markov chain with a  \emph{non-negative
convolution of projections }(Definition \ref{defn:nonnegativeconvolutionprojection}),
one can change the piece sizes in the breaking step to have any desired
distribution. For example, the ``remove and re-insert a box'' chain
above comes from the map $m(\Proj_{1}\otimes\iota)\Delta$, using
exactly the same Doob transform. ($\iota$ is the identity map and
$\Proj_{1}$ is projection to the subspace of degree 1.) On the shuffle
algebra, $m(\Proj_{1}\otimes\iota)\Delta$ defines the much studied
\emph{top-to-random shuffle}: take the top card off the deck, then
re-insert it at a uniformly chosen position. The other non-negative
convolutions of projections recover the shuffles of \parencite{cppriffleshuffle},
where the deck is cut in some specified distribution, and then the
cards dropped one by one from the bottom of piles chosen with probability
proportional to pile size. Their Corollaries 5.1 and 5.2 are a formula
for the composition of such shuffles, and an upper bound for the mixing
time. The inverses of such shuffles are examples of \emph{pop shuffles}
of \parencite{hyperplanewalk}, so the eigenvalues and multiplicities
can be calculated with their hyperplane walk theory. 

This new class of Markov chains admit analysis by the same techniques
as for Hopf-power Markov chains. Maps between Hopf algebras ``respecting
the bases'' again induce projections of their associated chains (Theorem
\ref{thm:projection} below). Consequently, the descent set is a Markov
statistic under all these shuffling schemes. Existing literature on
convolutions of projections provides the eigenvalues and multiplicities
of these transition matrices. In some cases, there are eigenbasis
algorithms resembling those for the Hopf-power chains. Once again,
the stationary distributions of all these chains are accessible -
they are precisely the same as those for the Hopf-power chains.

One notable shuffle outside this framework is \emph{random-to-random}:
uniformly choose a card to remove from the deck, and re-insert it
in a uniform position. Its defining linear map is an interesting operator
on other combinatorial Hopf algebras, and it would be great to find
a probability interpretation. 

This extended abstract is organised as follows: Section \ref{sec:prereqs}
gives the conditions on the two main characters in this story, the
state space basis of a combinatorial Hopf algebra and the non-negative
convolution of projections map. Section \ref{sec:chain} explains
how to construct the Markov chain via the Doob transform. Section
\ref{sec:theorems} states the three main theorems: how morphisms
of Hopf algebras lead to Markov statistics; the eigenvalues of the
transition matrices and their multiplicities; and the common stationary
distributions. Section \ref{sec:evectors} shows one scenario where
explicit eigenbasis formulae are available, and gives probability
applications both for shuffling and for a chain on trees.

Acknowledgements: I would like to thank Thomas Lam and Nathan Williams
for inspiring this project, and Nantel Bergeron, Persi Diaconis, Ian
Grojnowski, Philippe Nadeau and Franco Saliola for many helpful discussions.
The comments from the referees, especially the detailed pointers to
the literature, are much appreciated.

\section{Combinatorial Hopf Algebras and the Convolution Product\label{sec:prereqs}}

The starting point of our Markov chain construction is a combinatorial
Hopf algebra, which encodes the breaking and combining rules for our
family of combinatorial objects. An instructive example is the shuffle
algebra, whose associated Markov chains describe various models of
shuffling.
\begin{example}
\label{ex:shufflalg}The shuffle algebra has a basis $\calb$ of words,
which we will think of as decks of cards. For example, the word $accb$
will denote the deck with card $a$ on top, followed by two copies
of card $c$, and card $b$ on the bottom. This algebra is graded
by the lengths of the words, or the number of cards in the deck. The
product of two words is the sum of all their interleavings (with multiplicity),
and the coproduct of a word is the sum of all its deconcatenations.
For example, 
\[
m(ac\otimes cb)=2accb+acbc+cacb+cabc+cbac;
\]
\[
\Delta(accb)=1\otimes accb+a\otimes ccb+ac\otimes cb+acc\otimes b+accb\otimes1.
\]

\end{example}
The exposition of \parencite{vicreinernotes} gives background on combinatorial
Hopf algebras, and the opening of \parencite{foissychalist} contains an
extensive list of examples with references. Like many recent treatments,
these focus on generalisations of the symmetric functions, which,
though extremely important, are not so integral to the present Markov
chain application. The thesis of \parencite[Sec. 4]{combihopfalgintro}
is closer to the viewpoint herein. 

There is no rigorous definition of a combinatorial Hopf algebra. The
intuition is that such an algebra $\calh$ should have a basis $\calb$
indexed by a family of combinatorial objects, graded by their sizes.
(Assume throughout that the ground field is $\mathbb{Q}$ or $\mathbb{R}$,
to facilitate the probability applications.) Write $\calhn$ for the
subspace of $\calh$ of degree $n$, so $\calh=\bigoplus\calhn$.
Since the empty object is the only object of size 0, the vector space
$\calh$ is \emph{connected}, \textit{i.e.} $\dim\calh_{0}=1$. Now
equip $\calh$ with a multiplication map $m:\calh_{i}\otimes\calh_{j}\rightarrow\calh_{i+j}$:
for $w,z\in\calb$, set $m(w\otimes z)$ to be the sum (possibly weighted)
of all possible results from ``combining'' $w$ and $z$. Similarly,
the coproduct map $\Delta:\calhn\rightarrow\bigoplus_{i=0}^{n}\calh_{i}\otimes\calh_{n-i}$
takes $x\in\calb$ to the sum (possibly weighted) of $w\otimes z$
over all pairs $(w,z)$ obtainable by ``breaking'' $x$. There are
various axioms that these operations must satisfy. 

The above combinatorial interpretation of the product and coproduct
motivates that the structure constants be non-negative (conditions
i and ii below). This is crucial for the transition probabilities
in Definition \ref{defn:markovchain} to be non-negative. As we will
discover in Section \ref{sec:chain}, another important hypothesis
for the construction of the associated Markov chains is that $\Delta(x)\neq1\otimes x+x\otimes1$
for any $x\in\calb$ of degree greater than 1. Intuitively, it mandates
that ``every combinatorial object of size greater than 1 can be broken
into strictly smaller pieces''. Combining these criteria leads to
the following definition:
\begin{defn}
\label{defn:statespacebasis} \parencite[Def. 4.3.3]{mythesis} A basis
$\calb=\amalg_{n}\calbn$ of a graded connected Hopf algebra $\calh=\bigoplus_{n}\calhn$
is a \emph{state space basis} if:
\begin{enumerate}
\item for all $w,z\in\calb$, the expansion of $m(w\otimes z)$ in the $\calb$
basis has all coefficients non-negative;
\item for all $x\in\calb$, the expansion of $\Delta(x)$ in the $\calb\otimes\calb$
basis has all coefficients non-negative;
\item for $n>1$, the basis $\calbn$ contains no \emph{primitive} elements.
That is, $\Delta(x)\neq1\otimes x+x\otimes1$ for all $x\in\calbn$
with $n>1$.
\end{enumerate}
\end{defn}
Note that $\calh$ may contain primitive elements of high degree,
so long as they are not in the basis $\calb$. 

Having fixed a combinatorial Hopf algebra and a degree $n$, the next
step is to choose the distribution of sizes of pieces in the breaking
step of the Markov chain. As Step 1 of the interpretation in Section
\ref{sec:chain} will show, these distributions are in bijection,
up to scaling, with the class of non-negative convolutions of projections,
as defined below. For example, the operator $\Proj_{d_{1}}*\Proj_{d_{2}}$
will constrain the first piece to be of size exactly $d_{1}$, and
the second piece to be of size exactly $d_{2}$.

Given maps $f,g:\calh\rightarrow\calh$, their \emph{convolution product}
$f*g:\calh\rightarrow\calh$ is the composition 
\[
f*g:=m(f\otimes g)\Delta.
\]
Since the product $m$ on $\calh$ is associative and the coproduct
is coassociative, the convolution product is associative, and it is
useful to view 
\[
f_{1}*\dots*f_{a}\quad\mbox{as}\quad m^{[a]}(f_{1}\otimes\dots\otimes f_{a})\Delta^{[a]}.
\]
Here $m^{[a]}:\calh^{\otimes a}\rightarrow\calh$ and $\Delta^{[a]}:\calh\rightarrow\calh^{\otimes a}$
are the $a$-fold product and $a$-fold coproduct respectively, describing
the combining of $a$ objects and breaking one object into $a$ pieces
(some of which may be trivial). So $m^{[2]}=m$ and $\Delta^{[2]}=\Delta$.
(A precise definition, by induction, is as follows: $m^{[1]}:=\iota$,
$\proda:=m(m^{[a-1]}\otimes\iota)$; $\Delta^{[1]}:=\iota$, $\coproda:=(\iota\otimes\dots\otimes\iota\otimes\Delta)\Delta^{[a-1]}$.)

The present Markov chain application uses only the case $f_{i}=\Proj_{d_{i}}$,
the projection to the subspace of degree $d_{i}$. In other words,
$\Proj_{d}:\calh\rightarrow\calh_{d}$ is the linear map satisfying
$\Proj_{d}(x)=x$ for $x\in\calh_{d}$, and $\Proj_{d}(\calh_{i})=0$
if $i\neq d$. It will be convenient to allow the case $d=0$.

\parencite{descentoperators} studied linear combinations of these convolutions
of projections. He called them \emph{descent operators}, since his
Theorem II.7 shows that, on a commutative or cocommutative Hopf algebra,
these operators form an algebra under composition isomorphic to the
descent algebra of the symmetric group. This connection is central
to the eigenbasis algorithms for the present Markov chains, which
come from the following subset of these operators:
\begin{defn}
\label{defn:nonnegativeconvolutionprojection}Let $\calh$ be a graded
connected Hopf algebra, and fix an integer $a\geq2$. A map $\cpp:\calh\rightarrow\calh$
of the form 
\[
\cpp:=\sum_{d_{1},\dots,d_{a}}\alpha_{(d_{1},\dots,d_{a})}\Proj_{d_{1}}*\dots*\Proj_{d_{a}}
\]
is a \emph{non-negative convolution of projections on $\calhn$} if
\begin{enumerate}
\item for all \emph{weak-compositions} $D:=(d_{1},\dots,d_{a})$ of $n$
(that is, $d_{i}\geq0$ and $\sum d_{i}=n$), the coefficients $\alpha_{D}$
are non-negative; 
\item for some weak-composition $D$ of $n$ where each $d_{i}\neq n$,
the coefficient $\alpha_{D}$ is positive.
\end{enumerate}
\end{defn}
The second axiom ensures that $\cpp:\calhn\rightarrow\calhn$ is not
merely multiplication by a constant. Note that the map $\cpp$ does
not uniquely determine the coefficients $\alpha_{D}$, because of
the possibility of parts of size 0. Different choices of $\alpha_{D}$
lead to different interpretations of the same Markov chain. One final
remark on this definition: it is fine for infinitely many $\alpha_{D}$
to be non-zero, as the image $\cpp(x)$ of any particular $x\in\calh$
is always a finite sum. This is because there are only finitely many
weak-compositions of a given integer into exactly $a$ parts.

One key example of a non-negative convolution of projections is the
$a$th Hopf-power map $\proda\coproda$. This is the $a$th convolution
power of the identity map, so it corresponds to setting $\alpha_{D}\equiv1$
for all weak-compositions $D$. As the three-step interpretation of
Section \ref{sec:chain} will explain, the associated Markov chains
have a symmetric multinomial breaking step. This is the case previously
studied in \parencite{hopfpowerchains} and in \parencite{mythesis}. Another
important specialisation comes from $\alpha_{(1,n-1)}=1$, $\alpha_{D}=0$
if $D\neq(1,n-1)$ for any $n$, so $\cpp=\Proj_{1}*\iota$. This
map produces Markov chains which break off a singleton and reattach
it, analogous to the top-to-random shuffle of the introduction.

More examples are at the end of the next section.

\section{Building The Markov Chain\label{sec:chain}}

The following fact is the main motivation for the definition of a
Markov chain for each non-negative convolution of projections: the
probability that a riffle-shuffle takes a deck $x$ of $n$ cards
to a deck $y$ is the coefficient of $y$ in $\frac{1}{2^{n}}m\Delta(x)$.
In other words, the \emph{transition matrix} of the riffle-shuffling
of $n$ cards is $\left[\frac{1}{2^{n}}m\Delta\right]_{\calbn}^{T}$,
the transpose of the matrix of the linear operator $\frac{1}{2^{n}}m\Delta$
with respect to the basis $\calbn$ of words. A similar direct calculation
shows that the top-to-random shuffle of $n$ cards has transition
matrix $\left[\frac{1}{n}\Proj_{1}*\iota\right]_{\calbn}^{T}$. 

These observations suggest defining the transition matrix to be $\left[\frac{1}{\beta}\cpp\right]_{\calbn}^{T}$
on other combinatorial Hopf algebras, for other non-negative convolution
of projections $\cpp$ and some appropriate number $\beta$. However,
such a matrix represents transition probabilities only when each of
its rows sums to 1. In other words, this naive generalisation fails
if the rows of $\left[\cpp\right]_{\calbn}^{T}$ do not sum to the
same number. One of the major findings of \parencite[Th. 3.4]{hopfpowerchains}
is that, when $\cpp$ is a Hopf-power map $\cpp=\proda\coproda$,
and $\calb$ is a state space basis, then it is possible to define
a rescaling $\hatcalb$ of $\calb$ so the row sums of $\left[\cpp\right]_{\hatcalbn}^{T}$
are equal. \parencite[Sec. 4.3]{mythesis} gives a much slicker and more
general description of this rescaling, in terms of the Doob transform.
This allows a generalisation to linear maps on $\calh$ that are not
the Hopf-power. Indeed, \parencite[Th. 3.1]{mythesis} describes the choice
of rescalings $\hatcalb$ of $\calb$ that are available for arbitrary
linear maps.

It happens that the standard rescaling for the Hopf-power maps also
applies to non-negative convolutions of projections;  indeed, the
rescaling necessary to construct a $\cpp$-Markov chain depends only
on the underlying Hopf algebra, and not on the map $\cpp$. This resulting
Markov chain is:
\begin{defn}
\label{defn:markovchain}Let $\calh=\oplus_{n\geq0}\calhn$ be a graded
connected Hopf algebra with state space basis $\calb$. For $x\in\calbn$,
let $\eta(x)$ denote the sum of the coefficients (in the $\calb^{\otimes n}$
basis) of $\Proj_{1}^{\otimes n}\Delta^{[n]}(x)$, and let 
\[
\hatcalbn:=\left\{ \left.\frac{x}{\eta(x)}\right|x\in\calbn\right\} .
\]
Let $\cpp=\sum\alpha_{D}\Proj_{d_{1}}*\dots*\Proj_{d_{a}}$ be a non-negative
convolution of projections on $\calhn$. Then the\emph{ $\cpp$-Markov
chain on $\calbn$} has transition matrix
\[
\left[\frac{1}{\beta_{n}}\cpp\right]_{\hatcalbn}^{T},\quad\mbox{where}\quad\beta_{n}:=\sum_{d_{1}+\dots+d_{a}=n}\alpha_{D}\binom{n}{d_{1}\dots d_{a}}.
\]

\end{defn}
The hypotheses of a state space basis ensure that $\eta(x)$ is never
zero. In the shuffle algebra, $\eta(x)=1$ always, so no rescaling
is necessary.

Expressing the above transition probabilities in terms of a ``natural''
process requires careful analysis of the underlying Hopf algebra.
Fortunately, one only needs to do this once for each Hopf algebra
to interpret all its $\cpp$-Markov chains, as an analogue of \parencite[Th. 4.4.1]{mythesis}
shows that these have the form:
\begin{elabeling}{00.00.0000}
\item [{1.}] Choose a weak-composition $\left(d_{1},\dots,d_{a}\right)$
of $n$ with probability $\frac{1}{\beta_{n}}\alpha_{D}\binom{n}{d_{1}\dots d_{a}}$.
\item [{2.}] Choose a way to break into pieces of sizes $d_{1},\dots,d_{a}$.
\item [{3.}] Choose a way to combine these pieces.
\end{elabeling}
Here, the probabilities of the choices in steps 2 and 3 depend only
on the Hopf algebra, not on $\cpp$. (The exact expressions for these
probabilities are unsightly and not instructive; the interested reader
may consult \parencite[Th. 4.4.1]{mythesis}.) For example, a \emph{$\cpp$-shuffle}
(the $\cpp$-Markov chain on the shuffle algebra) of $n$ cards is
the following:
\begin{elabeling}{00.00.0000}
\item [{1,2.}] Cut the deck into piles of sizes $\left(d_{1},\dots,d_{a}\right)$
with probability $\frac{1}{\beta_{n}}\alpha_{D}\binom{n}{d_{1}\dots d_{a}}$.
\item [{3.}] Drop one-by-one the bottommost card from a pile chosen with
probability proportional to the current pile size. 
\end{elabeling}
Aside from the GSR riffle-shuffle, its $a$-handed generalisation,
and the top-to-random shuffle, here are some additional notable $\cpp$-shuffles:
\begin{example}
\label{ex:biasedshuffle}\parencite{biasedshufflesamisoundpersi} give
the mixing time for shuffles with biased cuts, when the shuffler prefers
to take more cards in one hand than the other. Here, the probability
of cutting $i$ cards off the top of a deck of $n$ cards is the asymmetric
binomial, $q^{i}(1-q)^{n-i}\binom{n}{i}$, for some parameter $q\in(0,1)$.
The associated non-negative convolution of projections is 
\[
\cpp=\sum_{i=0}^{n}q^{i}(1-q)^{n-i}\Proj_{i}*\Proj_{n-i}.
\]
This has an obvious $a$-handed generalisation, with $a-1$ parameters.
Setting all parameters to $\frac{1}{a}$ then recovers the $a$-handed
riffle-shuffle (even though the associated non-negative convolution
of projections is then $a^{-n}\iota^{*a}$ instead of $\iota^{*a}$,
as these Markov chains depend on $\cpp$ only up to scaling).
\end{example}

\begin{example}
\label{ex:topmtorandom-bothkinds}\parencite[Sec. 2  and Sec. 6, Ex. 2]{cppriffleshuffle}
discuss two notions of top-$m$-to-random shuffles: $\cpp=\Proj_{m}*\iota$
corresponds to cutting off $m$ cards and re-inserting them randomly,
keeping their relative order, whilst $\cpp=\Proj_{1}^{*m}*\iota$
cuts $m$ cards off and inserts them randomly in any order. For both
schemes, they show that $\frac{n}{m}\log n$ shuffles suffice to randomise
the deck. 
\end{example}

\begin{example}
\label{ex:toporbottomtorandom}Taking $\cpp=\Proj_{1}*\iota+\iota*\Proj_{1}$
produces a shuffle where the pile sizes are $(1,n-1)$ or $(n-1,1)$,
each with probability $\frac{1}{2}$. In other words, flip a fair
coin, and perform a top-to-random shuffle if the coin shows heads,
and a bottom-to-random shuffle if it shows tails. This is the (symmetric)
\emph{top-or-bottom-to-random shuffle} of \parencite[Sec. 6, Ex 4]{cppriffleshuffle}.
It is easy to introduce an asymmetry here: for $q\in[0,1]$, take
$\cpp=q\Proj_{1}*\iota+(1-q)\iota*\Proj_{1}$. Setting $q=1$ then
recovers the top-to-random shuffle. Theorem \ref{thm:evectors} below
exhibits an eigenbasis for this map on cocommutative Hopf algebras. 
\end{example}
The Markov chains coming from the above choices of $\cpp$, on other
Hopf-algebras, are also interesting. 
\begin{example}
\label{ex:repsofsn} The irreducible representations of the symmetric
groups form a basis of a Hopf algebra, with product being external
induction, and coproduct coming from restriction to Young subgroups.
It is straightforward to adapt \parencite[Ex. 4.4.3]{mythesis} to give
the following description for each step of a $\cpp$-Markov chain,
starting from a representation $x$ of $\sn$: 
\begin{elabeling}{00.00.0000}
\item [{1.}] Choose a Young subgroup $\mathfrak{S}_{d_{1}}\times\dots\times\mathfrak{S}_{d_{a}}$
of $\sn$ with probability $\frac{1}{\beta_{n}}\alpha_{D}\binom{n}{d_{1}\dots d_{a}}$.
\item [{2.}] Restrict the starting state $x$ to the chosen subgroup.
\item [{3.}] Induce it back up to $\sn$, then pick an irreducible constituent
with probability proportional to the dimension of its isotypic component.
\end{elabeling}

In particular, the $\Proj_{1}*\iota$-chain is restricting to $\mathfrak{S}_{n-1}$,
inducing back to $\sn$, then choosing an irreducible constituent
as in step 3. This chain previously appeared in the work of \parencite{downupchains},
where it generates central limit theorems for character ratios. 

\end{example}
For a more involved example, see Example \ref{ex:trees} regarding
the $(q\Proj_{1}*\iota+(1-q)\iota*\Proj_{1})$-Markov chain on trees.

\section{Projection Theorem and Eigenvalue Multiplicities\label{sec:theorems}}

As outlined in the introduction, viewing Markov chains in this Hopf-theoretic
framework gives two useful consequences. The first is the construction
of Markov statistics from maps between Hopf algebras: 
\begin{thm}
\label{thm:projection}Let $\calh$, $\barcalh$ be graded, connected
Hopf algebras with bases $\calb$, \textup{$\barcalb$} respectively.
Suppose in addition that $\calb$ is a state space basis. Let $\cpp$
be a non-negative convolution of projections. If $\theta:\calh\rightarrow\barcalh$
is a Hopf-morphism such that $\theta(\calbn)=\barcalb_{n}$ for all
$n$, then the $\cpp$-chain on $\barcalb_{n}$ is the projection
under $\theta$ of the $\cpp$-chain on $\calb_{n}$. 
\end{thm}
It follows that $\theta$ is a Markov statistic for the $\cpp$-chain
on $\calbn$ - this fact would be interesting by itself even if the
projected chain were not identified as the $\cpp$-chain on the target
Hopf algebra. 

The second profit of the Hopf formulation is the following expression
for all the eigenvalues and multiplicites of these Markov chains,
which shed some light on their long term behaviour. It comes from
combining \parencite[Th. II.7]{descentoperators}, \parencite[Prop. 3.12]{ncsym2},
\parencite[Prop. 3.10]{cppdiagonalisable} and the arguments of \parencite[Th. 3.21]{ncsym2}
and of \parencite{diagonalisingusinggrh}.
\begin{thm}
\label{thm:evalues}Work in the setup of Definition \ref{defn:markovchain}.
Given a partition $\lambda:=(\lambda_{1},\dots,\lambda_{l})$ and
a weak-composition $D=(d_{1},\dots,d_{a})$, let $\langle\lambda,D\rangle$
denote the number of set partitions $B_{1}|\dots|B_{a}$ of $\{1,2,\dots,l\}$
such that, for each $i\in\{1,\dots a\}$, we have $\sum_{j\in B_{i}}\lambda_{j}=d_{i}$.
(So $\langle\lambda,D\rangle$ is equal to the inner product $\langle p_{\lambda},h_{D}\rangle$
of symmetric functions, hence the notation.) Then the eigenvalues
of $\frac{1}{\beta_{n}}\cpp:\calhn\rightarrow\calhn$ are 
\[
\left\{ \left.\frac{\beta_{\lambda}}{\beta_{n}}:=\frac{1}{\beta_{n}}\sum_{D\vdash n}\alpha_{D}\langle\lambda,D\rangle\right|\lambda\vdash n\right\} ,
\]
and the multiplicity of the eigenvalue $\frac{\beta_{\lambda}}{\beta_{n}}$
is the coefficient of $x_{\lambda}:=x_{\lambda_{1}}\dots x_{\lambda_{l}}$
in the generating function $\prod_{i}(1-x_{i})^{-b_{i}}$, where $b_{i}$
satisfies 
\[
\sum_{n}\dim\calhn t^{n}=\prod_{i}(1-t^{i})^{-b_{i}}.
\]
Futhermore, $\cpp$ is diagonalisable if $\calh$ is commutative or
cocommutative.
\end{thm}
Note that $\beta_{(n)}$ agrees with the $\beta_{n}$ of Definition
\ref{defn:markovchain}, so this is not a point of confusion. Here's
how this formula specialises to some key examples:
\begin{example}
\label{ex:toptorandom-evalues}Let $\cpp=\Proj_{1}*\iota$, the top-to-random
map. Recall that this corresponds to $\alpha_{(1,n-1)}=1$, and all
other $\alpha_{D}=0$. So, for all partitions $\lambda$ of $n$,
it holds that $\beta_{\lambda}=\langle\lambda,(1,n-1)\rangle$, and
this is the number of parts of size 1 in $\lambda$, which can be
$0,1,\dots,n-2,$ or $n$. So the eigenvalues of a top-to-random chain
on any Hopf algebra are $\beta_{\lambda}/\beta_{(n)}=0,\frac{1}{n},\frac{2}{n},\dots,\frac{n-2}{n},1$. 

In the case of the shuffle algebra, for a deck of distinct cards,
\parencite[Th. 4.1]{cppriffleshuffle} show that the multiplicity of the
eigenvalue $\frac{j}{n}$ is the number of permutations with $n-j$
fixed points, and find projection operators for each eigenspace (on
the right). \parencite[Sec. 4.6]{oneminuse} produce an eigenbasis by
associating each permutation with $n-j$ fixed points to an eigenvector
of eigenvalue $\frac{j}{n}$. The present $\Proj_{1}*\iota$-chain
framework generalises this eigenbasis algorithm for decks with repeated
cards; see the remark after Proposition \ref{prop:descentspeaks}.
It follows from a multigraded refinement of Theorem \ref{thm:evalues}
above that, for any $\cpp$-shuffle of a deck of distinct cards, the
multiplicity of the eigenvalue $\beta_{\lambda}/\beta_{(n)}$ is the
number of permutations of cycle type $\lambda$.
\end{example}

\begin{example}
\label{ex:toporbottomtorandom-evalues}Take $\cpp=(q\Proj_{1}*\iota+(1-q)\iota*\Proj_{1})$,
the ``asymmetric top-or-bottom-to-random'' operator. Its eigenvalues
are 
\[
\frac{\beta_{\lambda}}{\beta_{n}}=\frac{q\langle\lambda,(1,n-1)\rangle+(1-q)\langle\lambda,(n-1,1)\rangle}{n}.
\]
Note that the definition of $\langle\lambda,D\rangle$ depends only
on the part sizes of the composition $D$ and not on their order,
so $\langle\lambda,(1,n-1)\rangle=\langle\lambda,(n-1,1)\rangle$.
Hence the eigenvalues of top-or-bottom-to-random are the same as for
top-to-random in Example \ref{ex:toptorandom-evalues} above, with
the same multiplicities.
\end{example}
Using the multiplicity in Theorem \ref{thm:evalues} for the eigenvalue
$1=\frac{\beta_{(n)}}{\beta_{n}}$ shows that the following expressions,
which are easily shown to be linearly independent stationary distributions,
do span the eigenspace of eigenvalue 1.
\begin{thm}
\label{thm:stationarydistribution}For a fixed state space basis $\calb$,
all $\cpp$-Markov chain on $\calbn$ have the same set of stationary
distributions. These can be uniquely written as a linear combination
of the functions 
\[
\pi_{c_{1},\dots,c_{n}}(x):=\frac{\eta(x)}{n!^{2}}\sum_{\sigma\in\sn}\mbox{coefficient of }x\mbox{ in the product }c_{\sigma(1)}\dots c_{\sigma(n)}
\]
for each multiset $\{c_{1},\dots,c_{n}\}$ in $\calb_{1}$. (Here,
$\eta(x)$ are the rescaling constants of Definition \ref{defn:markovchain}.)
\end{thm}
As noted in \parencite[Th. 4.5.1]{mythesis}, $\pi_{c_{1},\dots,c_{n}}(x)$
essentially enumerates the ways to build $x$ out of $c_{1},\dots,c_{n}$
using the multiplication of the combinatorial Hopf algebra, and to
then break it into singletons. In the case of card-shuffling, the
unique stationary distribution for all $\cpp$-shuffles is the uniform
distribution.

\section{Eigenvectors and Applications\label{sec:evectors}}

Since the coefficients $\alpha_{D}$ of a non-negative convolution
of projections can take any non-negative value, it's not surprising
that there is no neat universal eigenbasis algorithm for these chains.
However, one case which works out nicely, thanks to the theory of
dual graded graphs of \parencite[Th. 1.6.6]{dgg}, is the top-or-bottom-to-random
chain of Example \ref{ex:toporbottomtorandom}:
\begin{thm}
\label{thm:evectors}Let $\calh$ be a graded connected Hopf algebra,
and $\calp$ be a (graded) basis of its primitive subspace. Write
$\calp$ as the disjoint union $\calp_{1}\amalg\calp_{>1}$, where
$\calp_{1}$ has degree 1. Set 
\[
E_{j}:=\left\{ \sum_{i=0}^{j}\sum_{\sigma\in\sj}\binom{j}{i}q^{i}(1-q)^{j-i}c_{\sigma(1)}\dots c_{\sigma(i)}\left(\sum_{\tau\in\sk}p_{\tau(1)}\dots p_{\tau(k)}\right)c_{\sigma(i+1)}\dots c_{\sigma(j)}\right\} ,
\]
ranging over all multisets $\{c_{1},\dots,c_{j}\}$ of $\calp_{1}$,
and all multisets $\{p_{1},\dots,p_{k}\}$ of $\calp_{>1}$ where
$\deg p_{1}+\dots+\deg p_{k}=n-j$. Then $E_{j}$ is a linearly independent
set of eigenvectors of eigenvalue $\frac{j}{n}$ for the operator
$\frac{1}{n}(q\Proj_{1}*\iota+(1-q)\iota*\Proj_{1}):\calhn\rightarrow\calhn$.
Furthermore, if $\calh$ is cocommutative, then $\amalg_{j=0}^{n-2}E_{j}\amalg E_{n}$
is a basis.
\end{thm}
Here are some simple applications of these eigenvectors to the top-or-bottom-to-random
shuffle of a deck of distinct cards, analogous to the statements for
riffle-shuffling in the fifth paragraph of the introduction. (The
shuffle algebra is commutative, so its dual is cocommutative, and
the eigenvectors of $\cpp$ on this dual give right eigenfunctions
of the transition matrix, from which one deduces these results.)
\begin{prop}
\label{prop:descentspeaks}Let $\{X_{t}\}$ denote the top-or-bottom-to-random
shuffle, with parameter $q$, of a deck of $n$ distinct cards. Suppose
the starting deck $X_{0}$ is in ascending order. Let $\Des(X)\subseteq\{1,2,\dots,n-1\}$
and $\Peak(X)\subseteq\{1,2,\dots,n-2\}$ denote the positions of
the descents and peaks of $X$ respectively. Then 
\begin{eqnarray*}
\mbox{Expectation}\left\{ \sum_{i\in\Des(X_{t})}\binom{n-2}{i-1}q^{i-1}(1-q)^{n-1-i}\right\}  & = & \left(1-\left(\frac{n-2}{n}\right)^{t}\right)\frac{1}{2};\\
\mbox{Expectation}\left\{ \sum_{i\in\Peak(X_{t})}\binom{n-3}{i-1}q^{i-1}(1-q)^{n-2-i}\right\}  & = & \left(1-\left(\frac{n-3}{n}\right)^{t}\right)\frac{1}{3}.
\end{eqnarray*}

\end{prop}
Setting $q=1$ in Theorem \ref{thm:evectors} gives the eigenvectors
for the top-to-random operator $\Proj_{1}*\iota$ as 
\[
\sum_{\sigma\in\sj}c_{\sigma(1)}\dots c_{\sigma(j)}\sum_{\tau\in\sk}p_{\tau(1)}\dots p_{\tau(k)}.
\]
These are also the eigenvectors of $\Proj_{1}^{*m}*\iota$, the unordered
version of top-$m$-to-random as in Example \ref{ex:topmtorandom-bothkinds},
with eigenvalue $\binom{j}{m}/\binom{n}{m}$. The reason is that,
on a cocommutative Hopf algebra, the map $\Proj_{1}^{*m}*\iota$ is
a polynomial in $\Proj_{1}*\iota$. Similarly, the $E_{j}$ in Theorem
\ref{thm:evectors} are eigenvectors of any polynomial in $q\Proj_{1}*\iota+(1-q)\iota*\Proj_{1}$.
In particular, they have eigenvalue $q_{2}^{j}$ for the following
map, corresponding to the \emph{trinomial-top-and-bottom-to-random
shuffle} of \parencite[Sec. 6, Ex 6]{cppriffleshuffle}: 
\[
\sum_{m_{1}+m_{2}+m_{3}=n}\frac{1}{m_{1}!m_{3}!}q_{1}^{m_{1}}q_{2}^{m_{2}}q_{3}^{m_{3}}\Proj_{1}^{*m_{1}}*\iota*\Proj_{1}^{*m_{3}},
\]
(Here, $q_{1}+q_{2}+q_{3}=1$, and the previous $q$ is $\frac{q_{1}}{q_{1}+q_{3}}$
in terms of the new parameters.)

To finish, here is an example away from the world of card-shuffling,
to illustrate the diversity of Markov chains that this framework can
analyse.
\begin{example}
\label{ex:trees}We will study the trinomial-top-and-bottom-to-random
Markov chain (of two paragraphs prior) on the Connes-Kreimer Hopf
algebra of rooted forests. Take as state space basis the set of all
rooted forests - that is, a disjoint union of trees, each of which
has a distinguished root vertex. (The vertices are unlabelled, and
the embedding of the tree in the plane is immaterial.) The degree
of a forest is its number of vertices. The product of two trees is
their disjoint union, and the coproduct of a tree $T$ is the sum
of $T\backslash S\otimes S$ over all connected subtrees $S$ of $T$
which are either empty or contain the root of $T$. Hence $\iota*\Proj_{1}$
corresponds to removing a root of a forest, and $\Proj_{1}*\iota$
to removing a leaf. For a full definition of this algebra, see \parencite[second half of Sec. 2]{cktrees2}.

The Hopf-power Markov chain on this algebra was the subject of \parencite[Sec. 5.3]{mythesis}.
Adapting Theorem 5.3.8 there gives the following description of the
trinomial-top-and-bottom-to-random chain: 

Suppose a company has a forest structure, so all employees have at
most one direct superior. All employees are either regular employees
or VPs, and the superior of a VP is necessarily also a VP (so the
VPs in each connected component of the company form a subtree containing
the root). 

Each month, each regular employee independently produces excellent
work with probability $q_{3}$, average work with probability $q_{2}$,
and subpar work with probability $q_{1}$ (where $q_{1}+q_{2}+q_{3}=1$.)
For each employee who produced excellent work, one by one in a random
order, the regular employee furthest up the chain of superiority from
him becomes a VP. Then the employees who produced subpar work are
fired, one by one starting from the most superior. Each firing causes
a cascade of promotions: first, someone further down the chain of
superiority from the fired employee is uniformly selected to replace
him. Then, if the promoted employee was superior to anyone, then one
of those is uniformly selected and promoted to his position. This
process continues until someone who is not superior to anyone is promoted. 

The chain keeps track of the structure of the regular employees, but
does not know which employee is taking which position in the forest
structure, nor does it see the structure of the VPs. 

The cases $j=2,3$ below are analogues of Proposition \ref{prop:descentspeaks}
for this chain. These are inequalities, rather than equalities, because
the exact eigenvectors in Theorem \ref{thm:evectors} are very complicated
(involving a second sum), so to obtain a slicker result, we use instead
the estimates $\fo_{j}$. Defining these requires some more notation:
for a regular employee $u$, let $\desc(u)$ (resp. $\anc(u)$) denote
the number of regular employees who are further down (resp. up) from
$u$ in the chain of superiority, including $u$ himself in both counts.
(In tree language, these are the descendants and the ancestors). Also,
write $n'(u)$ for the size of the connected component of regular
employees containing $u$. 
\begin{prop}
\label{prop:tres-efns}Let $\{X_{t}\}$ denote the trinomial-top-and-bottom-to-random
Markov chain on the Connes-Kreimer Hopf algebra of rooted forests,
with interpretation and notation as above. For each integer $j\geq2$,
define the following functions on forests: 
\[
\fo_{j}(T):=\sum_{u\in T}q_{1}^{\desc(u)}q_{3}^{\anc(u)}\binom{\desc(u)}{j}.
\]
(The binomial coefficient is 0 if $\desc(u)<j$.) Then 
\[
\mbox{Expectation}\left\{ \fo_{j}(X_{t})\right\} \leq q_{2}^{jt}\fo_{j}(X_{0})\max_{u\in X_{0}:\desc(u)\geq j}\left\{ \binom{n'(u)}{\anc(u)-1}\right\} .
\]
 
\end{prop}

\end{example}
\printbibliography

\end{document}